\documentclass[english]{article}
\usepackage[T1]{fontenc}
\usepackage[latin9]{inputenc}
\usepackage{amsmath}
\usepackage{amsthm}
\usepackage{verbatim}
\usepackage{amssymb}

\makeatletter
\theoremstyle{plain}
    \ifx\thechapter\undefined
	    \newtheorem{thm}{\protect\theoremname}
	  \else
      \newtheorem{thm}{\protect\theoremname}[chapter]
    \fi
\theoremstyle{plain}
\newtheorem*{lem*}{\protect\lemmaname}
\theoremstyle{remark}
\newtheorem*{rem*}{\protect\remarkname}

\usepackage{geometry}

\makeatother

\usepackage{babel}
\providecommand{\lemmaname}{Lemma}
\providecommand{\remarkname}{Remark}
\providecommand{\theoremname}{Theorem}

\begin{document}
\title{  Note on the Coincidence theorem}
 \author{Rado\v s Baki\'{c}}
\date{}
\maketitle
\begin{abstract} 
We are proving Coincidence theorem due to Walsh for the case when total degree of a polynomial is less then number of arguments. As an application we prove generalizations of the classical composition theorem. Also, the following result has been proven: if $p(z)$ is a complex polynomial of degree $n$, then closed disk D that contains at least $n-1$ of its zeros (counting multiplicity) contains at least $\left[\frac{n-2k+1}{2} \right]$   zeros  of its k-th derivative, provided that arithmetical mean of these zeros is also center of D.  \end{abstract}

\vspace{1cm}

 Let $a(z)=\sum_{k=0}^{n}a_{k}z^{k}$ and $b(z)=\sum_{k=0}^{n}b_{k}z^{k}$
 be two complex polynomials of degree $n$. For them we can define linear operator $ A(a,b)=\sum_{k=0}^{n}(-1)^{k}\frac{  a_{k}b_{n-k}  }{ {n\choose k}  }$.
 
  If $A(a,b)=0$, then $a$ and $b$ are said to be apolar polynomials. For apolar polynomials holds classical theorem due to Grace:
  
  \textbf{Theorem of Grace:} If all zeros of a polynomial $a(z)$ are contained in some circular region $S$, then at least one zero of $b(z)$ is contained in $S$, provided that $a(z)$ and $b(z)$ are apolar.
 
 Circular region is (open or closed)  disk or half-plane or their exterior.   Some generalizations of the  theorem of Grace can be found in [4].

 One consequence of this theorem is a classical Coincidence theorem due to Walsh:
 
 \textbf{ Coincidence theorem:} Let $p(z_{1},z_{2},\ldots,z_{n})$ be a symmetric complex polynomial of total degree $n$, and of degree 1 in each $z_{i}$. Suppose that  $w_{1},w_{2},\ldots,w_{n}$ are complex nu,mbers that are contained in some circular region $S$.  Then exists $z\in S$ such that $p(w_{1},w_{2},\ldots,w_{n})=p(z,z,\ldots,z)$. 
 
 In his papers [1], [2] and [3] Aziz showed that above theorem holds if total degree is less then $n$, provided that mentioned circular region is convex. In general, Coincidence theorem is not true is $S$ is exterior of a disk and total degree is less then $n$. For example, if we set $p(z_{1},z_{2})=z_{1}+z_{2}$, and $S$ is closed exterior of the unit disk, then equation $0=p(-1,1)=2z$ has no solution in $S$. 
 
 We are now going to strengthen above   results   of Aziz and also extend it to the exterior of a disk. 
 
 \textbf{Theorem 1} Let $p(z_{1},z_{2},\ldots,z_{n})$ be a symmetric complex polynomial of total degree $m\le n$, and of degree 1 in all $z_{i}$. Suppose now that $w_{1},w_{2},\ldots,w_{n} $ are complex numbers, such that zeros of $ q(z)^{(n-m)}$ are contained in some circular region $S$, where $q(z)=\prod_{i=1}^{n}(z-w_{i})$. Then equation $p(w_{1},w_{2},\ldots,w_{n})=p(z,z,\ldots,z)$ has solution in $S$.
 
  \textbf{Proof:}
   We can assume that $p(w_1,w_2,\ldots,w_n)=0$. Otherwise we can replace $p(z_1,z_2,\ldots,z_n)$ with $p(z_1,z_2,\ldots,z_n)-p(w_1,w_2,\ldots,w_n).$ 
   By well known representation theorem, polynomial $p$ can be represented as a linear combination of $e_{k}$ where $e_{k}$ is a elementary symmetric polynomial in  $z_{1},z_{2},\ldots,z_{n}$ of degree $k(e_{0}=1)$. Hence,  $p(z_{1},z_{2},\ldots,z_{n})=\sum_{k=0}^{m}E_{k}e_{k} $, for some complex constants $E_{i}$. We can also write $q(z)=\prod_{k=1}^{n}(z-w_{k})=\sum_{k=0}^{n}(-1)^{n-k}e_{n-k}z^k$, and so  $q^{(i)}(z)= \sum_{k=0}^{n-i}(i+k)\cdots (1+k)(-1)^{n-k-i}e_{n-k-i}z^{k}$. It can be easily checked that for $m<n$   $$
 p(z_{1},z_{2},\ldots,z_{n})=\sum_{k=0}^{m}E_{k}e_{k}=\frac{1}{n(n-1)\cdots (m+1)  }\sum_{k=0}^{m}\frac{E_{k}{n\choose k}(n-k)\cdots(m-k+1)e_{k}  }{{m\choose k}  }.  
  $$
 Then,  $$ p(w_{1},w_{2},\ldots,w_{n})=\frac{1}{n(n-1)\cdots(m+1) }   \sum_{k=0}^{m}\frac{E_{k}{n\choose k}(n-k)\cdots(m-k+1)e_{k}  }{{m\choose k}  }=0$$ is equivalent to $A(q^{(n-m)},r)=0$ where $r(z)= \sum_{k=0}^{m}E_{k} {n\choose k}z^k= p(z ,z ,\ldots,z )$. That means that $q^{(n-m)}(z)$ and $p(z,z,\ldots,z)$ are apolar (this is  also true for $n=m$). Hence, equation $$ p(z ,z ,\ldots,z )=0$$ has solution in $S$, and the theorem is proved.
 
 Let us note that in case when $S$ is convex, Theorem 1 is stronger (in general)  than mentioned results of Aziz, since any convex set that contains  zeros of some polynomial contains also its critical points due to Gauss-Lucas theorem. 
 
 Our next result is application of the Theorem 1. It generalizes   result given in [5]. 
 
 \textbf{Theorem 2.} Let $p(z)$ be a complex polynomial of degree $n$. Suppose that some closed disk $D$ contains $n-1$ zeros of $p(z)$. such that the centre of D is also the arithmetic mean of these zeros. Then disk D contains at least $\left[\frac{n-2k+1}{2} \right]  $ zeros of the k-th derivative of $p(z)$, where $ \left[\,\,\, \right] $ denotes integer part.
 
  \textbf{Proof:}    Let $ z_{1},z_{2},\ldots,z_{n}$ be zeros of $p(z)$ such that $ z_{1},z_{2},\ldots,z_{n-1}$ are contained in a disk D, and let $c$ be the centre of D, $c=\frac{1}{n-1}\sum_{k=1}^{n-1}z_{k}$. We can assume that $z_{n}$ is outside D, otherwise our theorem follows immediately from Gauss-Lucas theorem. Due to suitable rotation and translation, we can also assume that $z_{n}=0$ and $c$ is real and positive.
  
  If we find k-th derivatiove of  $p(z)=\prod_{k=1}^{n}(z-z_{k})$ it will be a sum of products of the form $(z-z_{l_1})\cdots (z-z_{l_{n-k}})$. We will group in one sum those products in which $z_{n}=0$ occurs, and rest of them in the other sum, i.e.: $$p^{(k)}(z)=z\Sigma_{1} +\Sigma_{2} $$
  
  Polynomial $p^{(k)}(z) $  can be viewed  as a polynomial in $ z_{1},z_{2},\ldots,z_{n-1}$ of total degree $n-k$  that satisfies conditions of the Theorem 1. If we set  $p^{(k)}(z)=z\sum_{1} +\sum_{2}=0, $ for some $z\in \mathbf{C}$, then by Theorem 1 exists $y\in D$ such that 
  
  \begin{align*}&p^{(k)}(z,z_{1},z_{2},\cdots, z_{n-1})\\&=p^{(k)}(z,y,y,\cdots, y)\\&=k!{n-1 \choose k} z(z-y)^{n-k-1}+k!{n-1 \choose k-1}(z-y)^{n-k}\\&= 0   \end{align*}
  
  and this is equivalent to $(z-y)^{n-k-1}(z-\frac{k}{n}y)=0.$ Hence, it follows that $p^{(k)}(z)=0 $ implies that either $z=y$(i.e. $ z\in D$)  or $z=\frac{k}{n}y$ for some $y\in D$. Now, let $w_{1},w_{2},\ldots,w_{n-k}$ be all zeros of $p^{(k)}(z)$. We can arrange these points such that first $m$ of them are in D, while all other points are outside D. So, all $w_{i}$ for $i>m$, are of the form $w_{i}=\frac{k}{n}y_{i}$ for some $y_{i}\in D$. Arithmetical means of zeros of $p(z)$ and $p^{(k)}(z)$ are equal, so $\frac{1}{n}\sum_{i=0}^{n}z_{i}=\frac{1}{n-k}\sum_{i=0}^{n-k}w_{i},$ i.e. $\frac{(n-1)c}{n} =\frac{1}{n-k}\sum_{i=0}^{n-k}w_{i}=\frac{1}{n-k}\left( \sum_{i=0}^{m}w_{i}+\frac{k}{n}\sum_{i=m+1}^{n-k}y_{i}\right),$ for some  $y_{i}\in D$.  All $w_i,\, i\le m,$ and $y_i$ in the above equation lie in D. So their real parts are less then $2c$. Therefore, if take real parts of both sides in the above equation, we obtain the following inequality: 
  $$ \frac{(n-1)c}{n}<\frac{2c}{n-k}(m+\frac{k}{n}(n-k-m))$$ and this is equivalent to  $\frac{n-2k-1}{2}<m$ i.e. $\left[\frac{n-2k+1}{2} \right] \le m,$ and the proof is completed. 
  
  \vspace{1cm}

Rado\v s Baki\'c

Teacher Education Faculty, University of Belgrade, Belgrade, Serbia

email: bakicr@gmail.com

\end{document}